\documentclass{article}
\usepackage{amssymb}
\usepackage{amsfonts}
\usepackage{amsmath}

\begin{document}

\title{Quantitative estimates in approximation by Bernstein-Durrmeyer-Choquet operators with respect to monotone and submodular set functions}
\author{Sorin G. Gal \\
Department of Mathematics and Computer Science, \\
University of Oradea, \\
Universitatii Street No.1, 410087, Oradea, Romania\\
E-mail: \textit{galso@uoradea.ro}}
\date{}
\maketitle
\begin{abstract}
For the qualitative results of pointwise and uniform approximation obtained in \cite{Gal-Opris}, we present general quantitative estimates in terms of the modulus of continuity and in terms of a $K$-functional, for the generalized multivariate Bernstein-Durrmeyer operator $M_{n, \Gamma_{n, x}}$, written in terms of the Choquet integral with respect to a family of monotone and submodular set functions, $\Gamma_{n, x}$, on the standard $d$-dimensional simplex. When $\Gamma_{n, x}$ reduces to two elements, one a Choquet submodular set function and the other one a Borel measure, for suitable modified Bernstein-Durrmeyer operators, univariate $L^{p}$-approximations, $p\ge 1$, with estimates in terms of a $K$-functional are proved. In the particular cases when $d=1$ and the Choquet integral is taken with respect to some concrete possibility measures, the pointwise estimate in terms of the modulus of continuity is detailed. Some simple concrete examples of operators improving the classical error estimates are presented. Potential applications to practical methods dealing with data, like learning theory and regression models, also are mentioned.
\end{abstract}
\textbf{AMS 2000 Mathematics Subject Classification}: 41A36, 41A25, 28A12, 28A25.

\textbf{Keywords and phrases}: Bernstein-Durrmeyer-Choquet operator, possibility measure, monotone and submodular set function, Choquet integral, quantitative estimates, pointwise approximation, uniform approximation, $L^{p}$-approximation, modulus of continuity, $K$-functional.

\section{Introduction}

The approximation properties of the multivariate Bernstein-Durrmeyer linear operator defined with respect to a nonnegative, bounded Borel measure $\mu:{\cal{B}}_{S^{d}}\to \mathbb{R}_{+}$, by
$$M_{n, \mu}(f)(x)$$
\begin{equation}\label{form1}
=\sum_{|\alpha|=n}\frac{\int_{S^{d}}f(t)B_{\alpha}(t)d\mu(t)}{\int_{S^{d}}B_{\alpha}(t)d\mu(t)}\cdot B_{\alpha}(x):=\sum_{|\alpha|=n}c(\alpha, \mu)\cdot B_{\alpha}(x), \, x\in S^{d},\, n\in \mathbb{N},
\end{equation}
where ${\cal{B}}_{S^{d}}$ denotes the sigma algebra of all Borel measurable subsets in the power set ${\cal{P}}(S^{d})$ and $f$ is supposed to be $\mu$-integrable on the standard simplex
$$S^{d}=\{(x_{1}, ..., x_{d}) ; 0\le x_{1}, ..., x_{d}\le 1, \, 0\le x_{1}+ ... +x_{d}\le 1\},$$
were studied in, e.g., the recent papers \cite{BJ}, \cite{Berd1}, \cite{Berd2}, \cite{Berd3} and \cite{Li}.

Note that in (\ref{form1}), it is used the notation
$$B_{\alpha}(x)=\frac{n!}{\alpha_{0}! \cdot \alpha_{1}! \cdot ... \cdot \alpha_{d}!}(1-x_{1}-x_{2}- ... - x_{d})^{\alpha_{0}}\cdot x_{1}^{\alpha_{1}} \cdot ... \cdot x_{d}^{\alpha_{d}}$$
$$=:\frac{n!}{\alpha_{0}! \cdot \alpha_{1}! \cdot ... \cdot \alpha_{d}!}\cdot P_{\alpha}(x),$$
where $\alpha=(\alpha_{0}, \alpha_{1}, ..., \alpha_{d})$, $\alpha_{j}\in \mathbb{N}\bigcup\{0\}$, $j=0, ..., d$, $|\alpha|=\alpha_{0}+\alpha_{1}+ ... +\alpha_{d}=n$.

In the very recent paper \cite{Gal-Opris}, we have proved that the approximation results in the above mentioned papers remain valid for the more general case when
$\mu$ is a monotone, normalized and submodular set function on $S^{d}$ and the integrals used in (\ref{form1}) are the nonlinear Choquet integrals with respect to $\mu$.

The main goal of this paper is to obtain quantitative estimates in terms of the modulus of continuity and in terms of some $K$-functionals, for
the pointwise and uniform approximation obtained in \cite{Gal-Opris} and for the univariate $L^{p}$-approximation, $p\ge 1$, in the case of more general multivariate Bernstein-Durrmeyer polynomial operators defined by
\begin{equation}\label{form1-bis}
M_{n, \Gamma_{n, x}}(f)(x)=\sum_{|\alpha|=n}c(\alpha, \mu_{n, \alpha, x})\cdot B_{\alpha}(x), \, x\in S^{d},\, n\in \mathbb{N},
\end{equation}
where
$$c(\alpha, \mu_{n, \alpha, x})=\frac{(C)\int_{S^{d}}f(t)B_{\alpha}(t)d\mu_{n, \alpha, x}(t)}{(C)\int_{S^{d}}B_{\alpha}(t)d\mu_{n, \alpha, x}(t)}
=\frac{(C)\int_{S^{d}}f(t)P_{\alpha}(t)d\mu_{n, \alpha, x}(t)}{(C)\int_{S^{d}}P_{\alpha}(t)d\mu_{n, \alpha, x}(t)}$$
and for every $n\in \mathbb{N}$ and $x\in S^{d}$, $\Gamma_{n, x}=(\mu_{n, \alpha, x})_{|\alpha|=n}$ is a family of bounded, monotone, submodular and strictly positive set functions on ${\cal{B}}_{S^{d}}$.

Note that if $\Gamma_{n, x}$ reduces to one element (i.e. $\mu_{n, \alpha, x}=\mu$ for all $n$, $x$ and $\alpha$), then the operator given by (\ref{form1-bis}) reduces to the operator considered in \cite{Gal-Opris}.

The plan of the paper is as follows. Section 2 contains some preliminaries on possibility theory and on Choquet integral. In Section 3, general quantitative estimates in terms of the modulus of continuity and in terms of a $K$-functional for the pointwise and uniform approximation by the operators $M_{n, \Gamma_{n, x}}(f)(x)$ defined by (\ref{form1-bis}) are obtained. Also, when $\Gamma_{n, x}$ reduces to two elements, one a Choquet submodular set function and the other one a Borel measure, for suitable modified Bernstein-Durrmeyer-Choquet operators, univariate $L^{p}$-approximations, $p\ge 1$, with quantitative estimates in terms of a $K$-functional are presented. Finally, in Section 4, in the particular case when $d=1$ and the Choquet integrals are taken with respect to some concrete possibility measures, the pointwise estimate in terms of the modulus of continuity is detailed. Also, some concrete example of operators improving the classical error estimates are presented and potential applications to practical methods dealing with data are mentioned.

\section{Preliminaries}

Firstly, we present a few known concepts in possibility theory useful for the next considerations.
For details, see, e.g., \cite{DubPrad}.

{\bf Definition 2.1.} For the non-empty set $\Omega$, denote by ${\cal{P}}(\Omega)$ the family of all subsets of $\Omega$.

(i) A function $\lambda : \Omega \to [0, 1]$ with the property $\sup\{\lambda(s) ; s\in \Omega\}=1$, is called possibility distribution on $\Omega$.

(ii) A possibility measure is a set function $P:{\cal{P}}(\Omega)\to [0, 1]$, satisfying the axioms $P(\emptyset)=0$, $P(\Omega)=1$ and $P(\bigcup_{i\in I}A_{i})=\sup\{P(A_{i}) ; i\in I\}$ for all $A_{i}\subset \Omega$, and any $I$, an at most countable family of indices. Note that if $A, B\subset \Omega$, $A\subset B$, then the last property easily implies that $P(A)\le P(B)$ and that $P(A\bigcup B)\le P(A)+P(B)$.

Any possibility distribution $\lambda$ on $\Omega$, induces the possibility measure $P_{\lambda}:{\cal{P}}(\Omega)\to [0, 1]$, given by the formula $P_{\lambda}(A)=\sup \{\lambda(s) ; s\in A\}$, for all $A\subset \Omega$ (see, e.g., \cite{DubPrad}, Chapter 1).

Some known concepts and results concerning the Choquet integral can be summarized by the following.

{\bf Definition 2.2.} Suppose $\Omega\not=\emptyset$ and let ${\cal{C}}$ be a $\sigma$-algebra of subsets in $\Omega$.

(i) (see, e.g., \cite{WK1}, p. 63) The set function  $\mu:{\cal{C}}\to [0, +\infty]$ is called a monotone set function (or capacity) if $\mu(\emptyset)=0$ and $\mu(A)\le \mu(B)$ for all $A, B\in {\cal{C}}$, with $A\subset B$. Also,
$\mu$ is called submodular if
$$\mu(A\bigcup B)+\mu(A\bigcap B)\le \mu(A)+\mu(B), \mbox{  for all } A, B\in {\cal{C}}.$$
$\mu$ is called bounded if $\mu(\Omega)<+\infty$ and normalized if $\mu(\Omega)=1$.

(ii) (see, e.g., \cite{WK1}, p. 233, or \cite{Choquet}) If $\mu$ is a monotone set function on ${\cal{C}}$
and if $f:\Omega \to \mathbb{R}$ is ${\cal{C}}$-measurable (that is, for any Borel subset $B\subset \mathbb{R}$ it follows $f^{-1}(B)\in {\cal{C}}$), then for any $A\in {\cal{C}}$, the concept of Choquet integral is defined by
$$(C)\int_{A} f d\mu=\int_{0}^{+\infty}\mu\left (F_{\beta}(f)\bigcap A\right )d\beta+\int_{-\infty}^{0}\left [\mu\left (F_{\beta}(f)\bigcap A\right )-\mu(A)\right ]d \beta,$$
where we used the notation $F_{\beta}(f)=\{\omega\in \Omega; f(\omega)\ge \beta\}$.
Notice that if $f\ge 0$ on $A$, then in the above formula we get $\int_{-\infty}^{0}=0$.

The function $f$ will be called Choquet integrable on $A$ if $(C)\int_{A}f d\mu\in \mathbb{R}$.

In what follows, we list some known properties of the Choquet integral.

{\bf Remark 2.3.} If $\mu:{\cal{C}}\to [0, +\infty]$ is a monotone set function, then the following properties hold :

(i) For all $a\ge 0$ we have $(C)\int_{A}af d\mu = a\cdot (C)\int_{A}f d\mu$ (if $f\ge 0$ then see, e.g., \cite{WK1}, Theorem 11.2, (5), p. 228 and if $f$ is of arbitrary sign, then see, e.g., \cite{Denn}, p. 64, Proposition 5.1, (ii)).

(ii) For all $c\in \mathbb{R}$ and $f$ of arbitrary sign, we have (see, e.g., \cite{WK1}, pp. 232-233, or \cite{Denn}, p. 65) $(C)\int_{A}(f+c)d \mu = (C)\int_{A}f d\mu + c\cdot \mu(A)$.

If $\mu$ is submodular too, then for all $f, g$ of arbitrary sign and lower bounded, we have
(see, e.g., \cite{Denn}, p. 75, Theorem 6.3)
$$(C)\int_{A}(f + g) d\mu \le (C)\int_{A}f d\mu + (C)\int_{A}g d\mu.$$

(iii) If $f\le g$ on $A$ then
$(C)\int_{A}f d\mu \le (C)\int_{A}g d\mu$ (see, e.g., \cite{WK1}, p. 228, Theorem 11.2, (3) if $f, g\ge 0$ and p. 232 if $f, g$ are of arbitrary sign).

(iv) Let $f\ge 0$. If $A\subset B$ then $(C)\int_{A}f d \mu \le (C)\int_{B}f d\mu.$
In addition, if $\mu$ is finitely subadditive, then
$$(C)\int_{A\bigcup B}f d\mu \le (C)\int_{A}f d\mu + (C)\int_{B}f d\mu.$$

(v) It is immediate that $(C)\int_{A}1\cdot d\mu(t)=\mu(A)$.

(vi) The formula $\mu(A)=\gamma(M(A))$, where
$\gamma :[0, 1]\to [0, 1]$ is an increasing and concave function, with $\gamma(0)=0$, $\gamma(1)=1$ and
$M$ is a probability measure (or only finitely additive) on a $\sigma$-algebra on $\Omega$ (that is, $M(\emptyset)=0$, $M(\Omega)=1$ and $M$ is countably additive), gives simple examples of normalized, monotone and submodular set functions (see, e.g., \cite{Denn}, pp. 16-17, Example 2.1). For example, we can take $\gamma(t)=\sqrt{t}$.

If the above $\gamma$ function is increasing, concave and satisfies only $\gamma(0)=0$, then for any bounded  Borel measure $m$, $\mu(A)=\gamma(m(A))$ gives a simple example of bounded, monotone and submodular set function.

Note that any possibility measure $\mu$ is normalized, monotone and submodular. Indeed, the axiom  $\mu(A\bigcup B)=\max\{\mu(A), \mu(B)\}$ implies the monotonicity, while the property $\mu(A\bigcap B)\le \min\{\mu(A), \mu(B)\}$ implies the submodularity.

(vii) If $\mu$ is a countably additive bounded measure, then the Choquet integral $(C)\int_{A}f d\mu$ reduces to the usual Lebesgue type integral (see, e.g., \cite{Denn}, p. 62, or \cite{WK1}, p. 226).

\section{Quantitative estimates for general Bernstein-Durrmeyer-Choquet operators}

Recall that $\mu:{\cal{B}}_{S^{d}}\to [0, +\infty)$ is said strictly positive if for every open set $A\subset \mathbb{R}^{n}$ with $A\cap S^{d}\not=\emptyset$, we have $\mu(A\cap S^{d})>0$.

The support of $\mu$ is defined by
$$supp(\mu)=\{x\in S^{d} ; \mu(N_{x})>0 \mbox{ for every open neighborhood } N_{x}\in {\cal{B}}_{S^{d}} \mbox{ of } x\}.$$
Note that the strict positivity of $\mu$, evidently implies the condition $supp(\mu)\setminus \partial S^{d}\not=\emptyset$, which guarantees that $(C)\int_{S^{d}}B_{\alpha}(t)d\mu(t)>0$, for all $B_{\alpha}$.

Let us consider $C_{+}(S^{d})=\{f : S^{d}\to \mathbb{R}_{+} ; f \mbox{ is continuous on } S^{d}\}$, endowed with the norm $\|F\|_{C(S^{d})}=\sup\{|F(x)| ; x\in S^{d}\}$.

The first main result of this section consists in the following general quantitative estimates in pointwise and uniform approximation.

{\bf Theorem 3.1.} {\it For each fixed $n\in \mathbb{N}$ and $x\in S^{d}$, let $\Gamma_{n, x}=\{\mu_{n, \alpha, x}\}_{|\alpha|=n}$ be a family of bounded, monotone, submodular and strictly positive set functions on ${\cal{B}}_{S^{d}}$.

(i) For every $f\in C_{+}(S^{d})$, $x=(x_{1}, ..., x_{d})\in S^{d}$, $n\in \mathbb{N}$, we have
$$|M_{n, \Gamma_{n, x}}(f)(x)-f(x)|\le 2\omega_{1}(f ; M_{n, \Gamma_{n, x}}(\varphi_{x})(x))_{S^{d}},$$
where $M_{n, \Gamma_{n, x}}(f)(x)$ is given by (\ref{form1-bis}), $\|x\|=\sqrt{\sum_{i=1}^{d}x_{i}^{2}}$, $\varphi_{x}(t)=\|t-x\|$ and $\omega_{1}(f ; \delta)_{S^{d}}=\sup\{|f(t)-f(x)| ; t, x\in S^{d}, \|t-x\|\le \delta\}$.

(ii) Suppose that the family $\Gamma_{n, x}$ does not depend on $x$. Then, for any $f\in C_{+}(S^{d})$ and $n\in \mathbb{N}$, we get
$$\|M_{n, \Gamma_{n}}(f)-f\|_{C(S^{d})}\le 2 K\left (f ; \frac{\Delta_{n}}{2}\right ),$$
where $\Delta_{n}=\sum_{i=1}^{d}\|M_{n, \Gamma_{n}}(|\varphi_{e_{i}}-x_{i}{\bf 1}|)\|_{C(S^{d})}$,
$$K(f; t)=\inf_{g\in C^{1}_{+}(S^{d})}\{\|f-g\|_{C(S^{d})} + t \|\nabla g\|_{C(S^{d})}\},$$
$C^{1}_{+}(S^{d})$ is the subspace of all functions $g\in C_{+}(S^{d})$ with continuous partial derivatives $\partial_{i} g$, $i \in \{1, ..., d\}$ and
$$\|\nabla g\|_{C(S^{d})}=\max_{i=\{1, ..., d\}}\{\|\partial_{i} g\|_{C(S^{d})}\},$$
$\varphi_{e_{i}}(x)=x_{i}$, $i\in \{1, ..., d\}$, $x=(x_{1}, ..., x_{d})$, ${\bf 1}(x)=1$, for all $x\in S^{d}$.}

{\bf Proof.} (i) For $x\in S^{d}$, $n\in \mathbb{N}$ and $|\alpha|=n$ arbitrary fixed, let us consider $T_{n, \alpha, x}:C_{+}(S^{d})\to \mathbb{R}_{+}$ defined by
$$T_{n, \alpha, x}(f)=(C)\int_{S^{d}}f(t) P_{\alpha}(t)d \mu_{n, \alpha, x}(t), f\in C_{+}(S^{d}).$$

Based on the above Remark 2.3, (i), (ii), (iii) and reasoning exactly as in the proof of Lemma 3.1 in \cite{Gal-Opris}, we get
$|T_{n, \alpha, x}(f)-T_{n, \alpha, x}(g)|\le T_{n, \alpha, x}(|f-g|)$. Then, since $T_{n, \alpha, x}$ is positively homogeneous, sublinear and monotonically increasing, we immediately get that $M_{n, \Gamma_{n, x}}$ keeps the same properties and as a consequence it follows
\begin{eqnarray}\label{eq1}
|M_{n, \Gamma_{n, x}}(f)(x)-M_{n, \Gamma_{n, x}}(g)(x)|\le M_{n, \Gamma_{n, x}}(|f-g|)(x),
\end{eqnarray}
$M_{n, \Gamma_{n, x}}(\lambda f)=\lambda M_{n, \Gamma_{n, x}}(f)$, $M_{n, \Gamma_{n, x}}(f+g)\le M_{n, \Gamma_{n, x}}(f)+M_{n, \Gamma_{n, x}}(g)$ and that $f\le g$ on $S^{d}$ implies $M_{n, \Gamma_{n, x}}(f)\le M_{n, \Gamma_{n, x}}(g)$ on $S^{d}$,
for all $\lambda\ge 0$, $f, g\in C_{+}(S^{d})$, $n\in \mathbb{N}$, $|\alpha|=n$, $x\in S^{d}$.

Denoting $e_{0}(t)=1$ for all $t\in S^{d}$, since obviously $M_{n, \Gamma_{n, x}}(e_{0})(x)=1$ for all $x\in S^{d}$ and taking into account the properties in Remark 2.3, (i) and (\ref{eq1}), for any fixed $x$ we obtain
\begin{eqnarray}\label{eq2}
|M_{n, \Gamma_{n, x}}(f)(x)-f(x)|&=&|M_{n, \Gamma_{n, x}}(f(t))(x)-M_{n, \Gamma_{n, x}}(f(x))(x)|\nonumber \\
&\le& M_{n, \Gamma_{n, x}}(|f(t)-f(x)|)(x).
\end{eqnarray}
But taking into account the properties of the modulus of continuity, for all $t,x\in S^{d}$ and $\delta>0$, we get
\begin{eqnarray}\label{eq3}
|f(t)-f(x)|\leq\omega_{1}(f;\|t-x\|)_{S^{d}}\leq\left[ \frac{1}{\delta }\|t-x\|+1
\right] \omega_{1}(f;\delta)_{S^{d}}.
\end{eqnarray}
Now, from (\ref{eq2}) and applying $M_{n, \Gamma_{n, x}}$ to (\ref{eq3}), by the properties of $M_{n, \Gamma_{n, x}}$ mentioned after the inequality (\ref{eq1}), we immediately get
$$|M_{n, \Gamma_{n, x}}(f)(x)-f(x)|\leq \left[ \frac{1}{\delta}M_{n, \Gamma_{n, x}}(\varphi_{x})(x)+1\right] \omega_{1}(f;\delta)_{S^{d}}.$$
Choosing here $\delta=M_{n, \Gamma_{n, x}}(\varphi_{x})(x)$, we obtain the desired estimate.

(ii) Let $f, g\in C_{+}(S^{d})$. We have
$$f(x)-M_{n, \Gamma_{n}}(f)(x)$$
$$=f(x)-g(x)+M_{n, \Gamma_{n}}(g)(x)-M_{n, \Gamma_{n}}(f)(x)+g(x)-M_{n, \Gamma_{n}}(g)(x),$$
which, by using (\ref{eq1}) too, implies
$$|f(x)-M_{n, \Gamma_{n}}(f)(x)|$$
$$\le |f(x)-g(x)|+|M_{n, \Gamma_{n}}(g)(x)-M_{n, \Gamma_{n}}(f)(x)|+|g(x)-M_{n, \Gamma_{n}}(g)(x)|$$
$$\le |f(x)-g(x)|+M_{n, \Gamma_{n}}(|g - f|)(x)+|g(x)-M_{n, \Gamma_{n}}(g)(x)|$$
$$\le 2 \|f-g\|_{C(S^{d})}+|g(x)-M_{n, \Gamma_{n}}(g)(x)|.$$
By following the lines in the proof of Theorem 4.5 in \cite{BJ}, since from the lines after relation (\ref{eq1}) in the above point (i), the operator $M_{n, \Gamma_{n}}$ is monotone and subadditive, for all $g\in C^{1}_{+}(S^{d})$, $x\in S^{d}$,
we immediately get
$$|g(x)-M_{n, \Gamma_{n}}(g)(x)|$$
$$\le M_{n, \Gamma_{n}}(|g- g(x){\bf 1}|)(x)\le \|\nabla g\|_{C(S^{d})}\cdot M_{n, \Gamma_{n}}\left (\sum_{i=1}^{d}|\varphi_{e_{i}}-x_{i}{\bf{1}}|\right )(x)$$
$$\le \|\nabla g\|_{C(S^{d})}\cdot \sum_{i=1}^{d} M_{n, \Gamma_{n}}\left (|\varphi_{e_{i}}-x_{i}{\bf{1}}|\right )(x)
\le \|\nabla g\|_{C(S^{d})}\cdot \Delta_{n}.$$
Concluding, it follows
$$\|f-M_{n, \Gamma_{n}}(f)\|_{C(S^{d})}\le 2 \left [\|f-g\|_{C(S^{d})} + \frac{\Delta_{n}}{2} \|\nabla g\|_{C(S^{d})}\right ],$$
which immediately implies the required estimate in (ii). $\hfill \square$

{\bf Remark 3.2.} The positivity of function $f$ in Theorem 3.1, (i) and (ii), is necessary because of the positive homogeneity of the Choquet integral used in their proofs. However, if $f$ is of arbitrary sign and lower bounded on $S^{d}$ with $f(x)-m\ge 0$, for all $x\in S^{d}$, then the statement of Theorem 3.1, (i), (ii), can be restated for the slightly modified Bernstein-Durrmeyer operator defined by
$$M_{n, \Gamma_{n, x}}^{*}(f)(x)=M_{n, \Gamma_{n, x}}(f-m)(x)+m.$$
Indeed, in the case of Theorem 3.1, (i), this is immediate from $\omega_{1}(f-m;\delta)_{S^{d}}=\omega_{1}(f;\delta)_{S^{d}}$ and from
$M_{n,\Gamma_{n, x}}^{*}(f)(x)-f(x)=M_{n,\Gamma_{n, x}}(f-m)(x)-(f(x)-m)$.
Note that in the case of Theorem 3.1, (ii), since we may consider here that $m<0$, we immediately get the relations
$$K(f - m ; t)=\inf_{g\in C^{1}_{+}(S^{d})}\{\|f-(g+m)\|_{C(S^{d})}+t \|\nabla g\|_{C(S^{d})}\}$$
$$=\inf_{g\in C^{1}_{+}(S^{d})}\{\|f-(g+m)\|_{C(S^{d})}+t \|\nabla (g+m)\|_{C(S^{d})}\}$$
$$=\inf_{h\in C^{1}(S^{d}), \, h\ge m}\{\|f-h\|_{C(S^{d})}+t \|\nabla h\|_{C(S^{d})}\}.$$

In the particular case when the family $\Gamma_{n, x}$ does not depend on $x$ and $n$, it is natural to ask for quantitative estimates in the $L^{1}$-approximation of the Choquet integrable functions (not necessarily continuous).
If, for example, $\Gamma_{n, x}=\{\mu\}$, $d=1$ and $\mu$ is a bounded, monotone and submodular set function, then for the Bernstein-Durrmeyer-Choquet operators
$$D_{n, \mu}(f)(x)=\sum_{k=0}^{n}p_{n, k}(x)\cdot \frac{(C)\int_{0}^{1}f(t)p_{n, k}(t)d \mu(t)}{(C)\int_{0}^{1}p_{n, k}(t)d \mu(t)}, \, p_{n, k}(x)={n\choose k} x^{k}(1-x)^{n-k},$$
with $f\in L^{1}_{\mu}$ meaning $f$ is ${\cal{B}}_{[0, 1]}$ measurable and $\|f\|_{L^{1}_{\mu}}=(C)\int_{0}^{1}|f(t)|d\mu(t)<\infty$, we get
$$\|D_{n, \mu}(f)\|_{L^{1}_{\mu}} \le \sum_{k=0}^{n}(C)\int_{0}^{1}p_{n, k}(t)|f(t)|d\mu(t)\le n \cdot \|f\|_{L^{1}_{\mu}}, \, n\in \mathbb{N}.$$
This is due to the fact that $(C)\int f d\mu$ is is not, in general, additive as function of $f$ (it is only subadditive).

Therefore, quantitative estimates in $L^{p}_{\mu}$-approximation by Bernstein-Durrmeyer-Choquet operators, remains, in the general case, an open question.

However, in the particular case when the family of set functions $\Gamma_{n, x}$ reduces, for example, to two elements (one being a Choquet submodular set function $\mu$ and the other one a Borel measure $\delta$), for suitable defined Bernstein-Durrmeyer-Choquet operators, quantitative $L^{p}_{\mu}$-approximation results, $p\ge 1$  hold. For this purpose, let us make the following notations :
$$L^{p}_{\mu}=\{f:[0, 1]\to \mathbb{R} ; f \mbox{ is } {\cal{B}}_{[0, 1]} \mbox{measurable and} (C)\int_{0}^{1}|f(t)|^{p}d\mu(t)<+\infty\},$$
$$L^{p}_{\mu, +}=L^{p}_{\mu}\bigcap \{f:[0, 1]\to \mathbb{R}_{+}\},$$
$$\overline{K}\left (f ; t\right )_{L^{p}_{\mu, \delta}}=\inf_{g\in C^{1}_{+}([0, 1])}\{\|f-g\|_{L^{p}_{\mu}}+\|f-g\|_{L^{p}_{\delta}} + t \|g^{\prime}\|_{C([0, 1])}\}.$$
It is easy to see that if $\mu\le \delta$ then for all $f\in L^{1}_{\delta}$ and $t\ge 0$, we have $2K(f; t/2)_{L^{p}_{\mu}}\le \overline{K}\left (f ; t\right )_{L^{p}_{\mu, \delta}}\le 2 K(f; t)_{L^{p}_{\delta}}$, where $K$ is of usual form and with the infimum taken for $g\in C^{1}_{+}([0, 1])$.

For $p=1$, we have :

{\bf Theorem 3.3.} {\it Let $\mu$ be a bounded, monotone, submodular and strictly positive set function on ${\cal{B}}_{[0, 1]}$ and $\delta$ a bounded, strictly positive Borel measure on ${\cal{B}}_{[0, 1]}$, such that $\mu(A)\le \delta(A)$ for all $A\in {\cal{B}}_{[0, 1]}$. Then, denoting  $L^{1}_{\delta, +}\subset L^{1}_{\mu, +}$ and defining
the Bernstein-Durrmeyer-Choquet operators
$$\overline{D}_{n, \delta, \mu}(f)(x)=\sum_{k=0}^{n-1}p_{n, k}(x)\cdot \frac{\int_{0}^{1}f(t)p_{n, k}(t)d \delta(t)}{\int_{0}^{1}p_{n, k}(t)d \delta(t)}+x^{n}\cdot \frac{(C)\int_{0}^{1}f(t)t^{n}d \mu(t)}{(C)\int_{0}^{1}t^{n}d \mu(t)},$$
for all $f\in L^{1}_{\delta, +}$, $n\in \mathbb{N}$ and denoting $\varphi_{x}(t)=|t-x|$, we have
$$\|f-\overline{D}_{n, \delta, \mu}(f)\|_{L^{1}_{\mu}}\le  2 \overline{K}\left (f ; \frac{\|\overline{D}_{n, \delta, \mu}(\varphi_{x})\|_{L^{1}_{\mu}}}{2}\right )_{L^{1}_{\mu, \delta}}.$$}
{\bf Proof.} Firstly, note that $\delta$ is monotone, submodular (in fact modular, i.e. submodular with equality)
strictly positive and that for all $f\in L^{1}_{\delta, +}$ we have $\int_{0}^{1}f(t)d \delta (t)=(C)\int_{0}^{1}f(t)d \delta(t)$ (see Remark 2.3, (vii)). From here, from $\mu\le \delta$ and from Definition 2.2, (ii), we immediately get $(C)\int_{0}^{1}f(t)d \mu(t) \le \int_{0}^{1}f(t)d \delta(t)$, which means $L^{1}_{\delta, +}\subset L^{1}_{\mu, +}$ and for all $f\in L^{1}_{\delta, +}$ implies
$$\|\overline{D}_{n, \delta, \mu}(f)\|_{L^{1}_{\mu}}\le \sum_{k=0}^{n-1}\frac{(C)\int_{0}^{1}p_{n, k}(x)d\mu(x)}{\int_{0}^{1}p_{n, k}(t)d\delta(t)}\cdot \int_{0}^{1}f(t)p_{n, k}(t)d \delta+(C)\int_{0}^{1}f(t)t^{n}d \mu(t)$$
\begin{equation}\label{eq111}
\le \int_{0}^{1}f(t)\left [\sum_{k=0}^{n-1}p_{n, k}(t)\right ]d \delta +(C)\int_{0}^{1}f(t)t^{n}d \mu(t)\le \|f\|_{L^{1}_{\delta}}+\|f\|_{L^{1}_{\mu}}.
\end{equation}
Let $f, g\in L^{1}_{\delta, +}$. From
$$|f(t)-\overline{D}_{n, \delta, \mu}(f)(t)|$$
$$\le |f(t)-g(t)|+|\overline{D}_{n, \delta, \mu}(g)(t)-\overline{D}_{n, \delta, \mu}(f)(t)|+|g(t)-\overline{D}_{n, \delta, \mu}(g)(t)|,$$
integrating with respect to $\mu$, from the properties of the Choquet integral, of the operator $\overline{D}_{n, \delta, \mu}$ (similar with those of $M_{n, \Gamma_{n, x}}$ in the proof of Theorem 3.1, (i)) and from (\ref{eq111}), we obtain
$$\|f-\overline{D}_{n, \delta, \mu}\|_{L^{1}_{\mu}}=(C)\int_{0}^{1}|f(t)-\overline{D}_{n, \delta, \mu}(f)(t)|d \mu(t)$$
$$\le (C)\int_{0}^{1}|f(t)-g(t)|d \mu(t) + (C)\int_{0}^{1}|\overline{D}_{n, \delta, \mu}(g)(t)-\overline{D}_{n, , \delta, \mu}(f)(t)|d \mu(t)$$
$$+(C)\int_{0}^{1}|g(t)-\overline{D}_{n, \delta, \mu}(g)(t)| d \mu(t)$$
$$\le \|f-g\|_{L^{1}_{\mu}} + (C)\int_{0}^{1} \overline{D}_{n, \delta, \mu}(|g-f|)(t)d \mu(t)+\|g - \overline{D}_{n, \delta, \mu}(g)\|_{L^{1}_{\mu}}$$
$$\le \|f-g\|_{L^{1}_{\mu}}+(\|f-g\|_{L^{1}_{\mu}}+\|f-g\|_{L^{1}_{\delta}})+\|g - \overline{D}_{n, \delta, \mu}(g)\|_{L^{1}_{\mu}}.$$
It remains to estimate $\|g-D_{n, \mu}(g)\|_{L^{1}_{\mu}}$.
But from
$$|g(x)-\overline{D}_{n, \mu}(g)(x)|=|\overline{D}_{n, \mu}(g(x))(x)-\overline{D}_{n, \mu}(g(t))(x)|\le \overline{D}_{n, \mu}(|g(x)-g(t)|)(x)$$
and since for $g\in C^{1}_{+}([0, 1])$, we get $|g(x)-g(t)|\le \|g^{\prime}\|_{C([0, 1])}|x-t|$, applying $\overline{D}_{n, \delta, \mu}$, it follows
$\overline{D}_{n, \delta, \mu}(|g(x)-g(t)|)(x)\le \|g^{\prime}\|_{C([0, 1])}\cdot \overline{D}_{n, \mu}(\varphi_{x})(x)$.

Therefore, integrating above with respect to $x$ and $\mu$, we obtain
$$\|g-\overline{D}_{n, \delta, \mu}(g)\|_{L^{1}_{\mu}}\le \|g^{\prime}\|_{C([0, 1])}\cdot \|\overline{D}_{n, \delta, \mu}(\varphi_{x})\|_{L^{1}_{\mu}},$$
which immediately leads to
$$\|f-\overline{D}_{n, \delta, \mu}(f)\|_{L^{1}_{\mu}}\le 2\|f-g\|_{L^{1}_{\mu}}+\|f-g\|_{L^{\delta}}+\|g^{\prime}\|_{C([0, 1])}\cdot \|\overline{D}_{n, \delta, \mu}(\varphi_{x})\|_{L^{1}_{\mu}}$$
$$\le 2\left (\|f-g\|_{L^{1}_{\mu}}+\|f-g\|_{L^{1}_{\delta}}+\|g^{\prime}\|_{C([0, 1])}\cdot \frac{\|\overline{D}_{n, \delta, \mu}(\varphi_{x})\|_{L^{1}_{\mu}}}{2}\right )$$
and to the conclusion of the theorem. $\hfill \square$

In what follows, because of some difference with respect to the case $p=1$, we extend separately  Theorem 3.3 to the $L^{p}_{\mu}$ space with $p>1$.

{\bf Theorem 3.4.} {\it Let $\mu$ be a bounded, monotone, submodular, strictly positive set function on ${\cal{B}}_{[0, 1]}$, which also is continuous by increasing sequences of sets, that is if $A_{n}\in {\cal{B}}_{[0, 1]}$, $n\in \mathbb{N}$, with $A_{n}\subset A_{n+1}$, for all $n$ and $A:=\bigcup _{n=1}^{\infty}A_{n}\in {\cal{B}}_{[0, 1]}$, then $\lim_{n\to \infty}\mu(A_{n})=\mu(A)$.

Also, let $\delta$ be a bounded, strictly positive Borel measure on ${\cal{B}}_{[0, 1]}$, such that $\mu(A)\le \delta(A)$ for all $A\in {\cal{B}}_{[0, 1]}$. Then, for any $p>1$,  $L^{p}_{\delta, +}\subset L^{p}_{\mu, +}$ and for
the Bernstein-Durrmeyer-Choquet operators $\overline{D}_{n, \delta, \mu}(f)(x)$ defined by Theorem 3.3,
for all $f\in L^{p}_{\delta, +}=L^{p}_{\delta}\bigcap \{f:[0, 1]\to [0, +\infty)\}$ and $n\in \mathbb{N}$, we have
$$\|f-\overline{D}_{n, \delta, \mu}(f)\|_{L^{p}_{\mu}}\le  2 \overline{K}\left (f ; \frac{\|\overline{D}_{n, \delta, \mu}(\varphi_{x})\|_{L^{p}_{\mu}}}{2}\right )_{L^{p}_{\mu, \delta}}.$$}
{\bf Proof.} The proof for $L^{p}_{\delta, +}\subset L^{p}_{\mu, +}$ follows exactly as in the proof of Theorem 3.3.

By the convexity of $t^{p}$ on $[0, +\infty)$, by $\sum_{k=0}^{n}p_{n, k}(x)=1$, we easily arrive at the inequalities (exactly as, for example, in the proof of Lemma 2.2 in \cite{Li})
$$\|\overline{D}_{n, \delta, \mu}(f)\|^{p}_{L^{p}_{\mu}}\le (C)\int_{0}^{1}\left [\sum_{k=1}^{n-1}p_{n, k}(x)\cdot
\frac{\left (\int_{0}^{1}f(t)p_{n, k}(t)d \delta(t)\right )^{p}}{\left (\int_{0}^{1}p_{n, k}(t)d \delta(t)\right )^{p}}\right .$$
$$\left . + x^{n}\cdot \frac{\left ((C)\int_{0}^{1}f(t)t^{n}d \mu(t)\right )^{p}}{\left ((C)\int_{0}^{1}t^{n}d \mu(t)\right )^{p}}\right ]d \mu(x).$$
Applying the  H\"older's inequality for the integrals from the denominators (in the case of Choquet integrals with respect to $\mu$, the inequality is the same with that for the integrals with respect to the Borel measure $\delta$, see. e.g., Theorem 3.5 in \cite{RSWang} or Theorem 2 in \cite{Cerda1})
and reasoning as in the proof of Lemma 2.2 in \cite{Li} and as for formula (\ref{eq111}) in the proof of Theorem 3.3,
we easily arrive at
$$\|\overline{D}_{n, \delta, \mu}(f)\|^{p}_{L^{p}_{\mu}}\le \|f\|_{L^{p}_{\mu}}+\|f\|_{L^{p}_{\delta}}, \mbox{ for all } f\in L^{p}_{\delta, +}.$$
Then, since H\"older's inequality for the Choquet integral with respect to $\mu$ implies the Minkowski inequality
(see, e.g., Theorem 3.7 in \cite{RSWang} or Theorem 2 in \cite{Cerda1}), using the above inequality and exactly the reasonings in the proof of Theorem 2.1 in \cite{Li}, we arrive at the desired inequality in the statement.

It remains to discuss the requirement on $\mu$ to be continuous by increasing sequences of sets. This is due
to the fact that for Choquet integrals, the H\"older's inequality hold only if both integrals from its right-hand side are not equal to zero (see the proofs of Theorem 3.5 in \cite{RSWang} or of Theorem 2 in \cite{Cerda1}).

To have valid the H\"older's inequality in its full generality, we need that for $F\ge 0$, $(C)\int_{0}^{1}F(t) d\mu=0$ if and only if $F(t)=0$, $\mu$ almost everywhere on $[0,1]$. But according to Theorem 11.3, p. 228 in \cite{WK1},
if $\mu$ is continuous by increasing sequences of sets, then the above mentioned property holds. $\hfill \square$

{\bf Remark 3.5.} Concrete choices for $\mu$ and $\delta$ in Theorem 3.3 can be, for example, $\delta(A)=m(A)$ and $\mu(A)=\sin[m(A)]$, where $m$ is the Lebesgue measure on ${\cal{B}}_{[0, 1]}$. Indeed, $\mu(A)\le m(A)$ for all $A\in {\cal{B}}_{[0, 1]}$ and since $\sin$ is concave on $[0,1]$ and $sin(0)=0$, by Remark 2.3, (vi) it follows that $\mu$ is bounded, monotone and submodular.

{\bf Remark 3.6.} It is easy to see that another Bernstein-Durrmeyer-Choquet operator satisfying the estimates in Theorems 3.3 and 3.4, can be defined by
$$\tilde{D}_{n, \delta, \mu}(f)(x)=(1-x)^{n}\frac{(C)\int_{0}^{1}f(t)(1-t)^{n}d \mu(t)}{(C)\int_{0}^{1}(1-t)^{n}d \mu(t)}+\sum_{k=1}^{n}p_{n, k}(x)\frac{\int_{0}^{1}f(t)p_{n, k}(t)d \delta(t)}{\int_{0}^{1}p_{n, k}(t)d \delta(t)}.$$
Also, defining
$$D^{\star}_{n, \delta, \mu}(f)(x)=(1-x)^{n}\frac{(C)\int_{0}^{1}f(t)(1-t)^{n}d \mu(t)}{(C)\int_{0}^{1}(1-t)^{n}d \mu(t)}+\sum_{k=1}^{n-1}p_{n, k}(x)\frac{\int_{0}^{1}f(t)p_{n, k}(t)d \delta(t)}{\int_{0}^{1}p_{n, k}(t)d \delta(t)}$$
$$+x^{n}\cdot \frac{(C)\int_{0}^{1}f(t)t^{n}d \mu(t)}{(C)\int_{0}^{1}t^{n}d \mu(t)},$$
by similar reasonings with those in the proofs of Theorems 3.3 and 3.4, for any $p\ge 1$ we immediately obtain the estimate
$$\|f-D^{\star}_{n, \delta, \mu}(f)\|_{L^{p}_{\mu}}\le 3 \overline{K}\left (f ; \frac{\|D^{\star}_{n, \delta, \mu}(\varphi_{x})\|_{L^{p}_{\mu}}}{3}\right ).$$

{\bf Remark 3.7.} For $\delta$ a bounded Borel measure on ${\cal{B}}_{[0, 1]}$, denote by $D_{n, \delta}$ the classical Bernstein-Durrmeyer operator (i.e. with all the integrals in terms of $\delta$). By Theorem 4.5 in \cite{BJ},
we have the estimate
$$\|D_{n, \delta}(f)-f\|_{L^{1}_{\delta}}\le 2K\left (f ; \frac{\|D_{n, \delta}(\varphi_{x})\|_{L^{1}_{\delta}}}{2}\right )_{L^{1}_{\delta}}, \mbox{ for all } f\in L^{1}_{\delta},$$
where $K(f ; t)_{L^{1}_{\delta}}=\inf_{g\in C^{1}([0, 1])}\{\|f-g\|_{L^{1}_{\delta}} + t \|g^{\prime}\|_{C([0, 1])}\}$.

Comparing with the estimate for $\|\overline{D}_{n, \delta, \mu}(f)-f\|_{L^{1}_{\mu}}$ in Theorem 3.3 and taking into account that $\|f-g\|_{L^{1}_{\mu}}\le \|f-g\|_{L^{1}_{\delta}}$ and $\|D_{n, \delta}(\varphi_{x})\|_{L^{1}_{\mu}}\le
\|D_{n, \delta}(\varphi_{x})\|_{L^{1}_{\delta}}$, it follows that it is possible that in some cases, $\overline{D}_{n, \delta, \mu}(f)$ in Theorem 3.3, approximates better $f\in L^{1}_{\delta, +}$ in the $L^{1}_{\mu}$-"norm" than approximate $D_{n, \delta}(f)$ the same function $f$ but in the $L^{1}_{\delta}$-norm.

\section{Concrete Bernstein-Durrmeyer-Choquet operators}

Since the estimates in Theorem 3.1 are of very general and abstract form, involving the apparently difficult to be calculated Choquet integrals, it is of interest to obtain in some particular cases, concrete expressions for the order of approximation.

In this sense, we will apply Theorem 3.1, (i), for $d=1$ and for some special choices of the submodular set functions.

Thus, we will consider the case of the measures of possibility. Denoting $p_{n, k}(x)={n \choose k}x^{k}(1-x)^{n-k}$, let us define $\lambda_{n, k}(t)=\frac{p_{n, k}(t)}{k^{k}n^{-n}(n-k)^{n-k}{n\choose k}}=\frac{t^{k}(1-t)^{n-k}}{k^{k}n^{-n}(n-k)^{n-k}}$, $k=0, ..., n$. Here, by convention we consider $0^{0}=1$, so that the cases $k=0$ and $k=n$ have sense.

By considering the root $\frac{k}{n}$ of $p^{\prime}_{n, k}(x)$, it is easy to see that $\max\{p_{n, k}(t) ; t\in [0, 1]\}=k^{k}n^{-n}(n-k)^{n-k}{n\choose k}$, which implies that each $\lambda_{n, k}$ is a possibility distribution on $[0, 1]$. Denoting by $P_{\lambda_{n, k}}$ the possibility measure induced by $\lambda_{n, k}$ and $\Gamma_{n, x}:=\Gamma_{n}=\{P_{\lambda_{n, k}}\}_{k=0}^{n}$ (i.e. $\Gamma$ is independent of $x$), the nonlinear Bernstein-Durrmeyer polynomial operators given by (\ref{form1-bis}), in terms of the Choquet integrals with respect to the set functions in $\Gamma_{n}$, will become
\begin{eqnarray}\label{eq8}
D_{n, \Gamma_{n}}(f)(x)=\sum_{k=0}^{n}p_{n, k}(x)\cdot \frac{(C)\int_{0}^{1}f(t)t^{k}(1-t)^{n-k}d P_{\lambda_{n, k}}(t)}{(C)\int_{0}^{1}t^{k}(1-t)^{n-k}d P_{\lambda_{n, k}}(t)}.
\end{eqnarray}
It is easy to see that any possibility measure $P_{\lambda_{n, k}}$ is bounded, monotone, submodular and strictly positive, $n\in \mathbb{N}$, $k=0, 1, ..., n$, so that we are under the hypothesis of Theorem 3.1, (i).

We can state the following result.

{\bf Theorem 4.1.} {\it If $D_{n, \Gamma_{n}}(f)(x)$ is given by (\ref{eq8}), then for every $f\in C_{+}([0, 1])$, $x\in [0, 1]$ and $n\in \mathbb{N}$, $n\ge 2$, we have
$$|D_{n, \Gamma_{n}}(f)(x)-f(x)|\le 2\omega_{1}\left (f ; \frac{(1+\sqrt{2})\sqrt{x(1-x)}+\sqrt{2}\sqrt{x}}{\sqrt{n}}+\frac{1}{n}\right )_{[0, 1]}.$$}
For its proof, we need the following auxiliary result.

{\bf Lemma 4.2.} {\it Let $n\in \mathbb{N}$, $n\ge 2$ and $x\in [0, 1]$. Denoting
$$A_{n, k}(x):=\sup\{|t-x|t^{k}(1-t)^{n-k};t\in [0, 1]\}=$$
$$\max\{\sup\{(t-x)t^{k}(1-t)^{n-k};t\in [x, 1]\}, \, \sup\{(x-t)t^{k}(1-t)^{n-k};t\in [0, x]\}\},$$
with the convention $0^{0}=1$, for all $k=0, ..., n$ we have
$$A_{n, k}(x)=\max\{(t_{2}-x)t_{2}^{k}(1-t_{2})^{n-k}, \, (x-t_{1})t_{1}^{k}(1-t_{1})^{n-k}\},$$
with $t_{1}, t_{2}$ given by
\begin{eqnarray}\label{eq5}
t_{1}=\frac{nx+k+1-\sqrt{\Delta}}{2(n+1)}, \, \, t_{2}=\frac{nx+k+1+\sqrt{\Delta}}{2(n+1)},
\end{eqnarray}
where
$$\Delta=(nx+k+1)^{2}-4k x(n+1)=n^{2}\left [(x+(k+1)/n)^{2}-4x\frac{k}{n}\cdot \frac{n+1}{n}\right ]$$
$$=(nx-k)^{2}+2x(n-k)+2k(1-x)+1\ge 1.$$}
{\bf Proof.} Let us denote $H_{n, k}(t)=t^{k}(1-t)^{n-k}|t-x|$, with $k\in \{0, ..., n\}$. We have two cases :
(i) $1\le k\le n-1$ and (ii) $k=0$ or $k=n$.

Case (i). For $t\in [x, 1]$ we obtain
$H_{n, k}(t)=(t-x)t^{k}(1-t)^{n-k}$ and from $H^{\prime}_{n, k}(t)=t^{k-1}(1-t)^{n-k-1}[-t^{2}(n+1)+t(nx+k+1)-kx]=0$,
it follows $-t^{2}(n+1)+t(nx+k+1)-kx=0$, which has the discriminant
$$\Delta=(nx+k+1)^{2}-4k x(n+1)=(nx-k)^{2}+2x(n-k)+2k(1-x)+1\ge 1.$$
Therefore, the quadratic equation has two real distinct solutions $t_{1}<t_{2}$
$$t_{1}=\frac{nx+k+1-\sqrt{\Delta}}{2(n+1)}, \, \, t_{2}=\frac{nx+k+1+\sqrt{\Delta}}{2(n+1)},$$
where by simple calculation we derive $0\le t_{1}<t_{2}\le 1$. Also, since $H_{n, k}(0)=H_{n, k}(x)=H_{n, k}(1)=0$ and $H_{n, k}(t)\ge 0$ for $t\in [x, 1]$, simple graphical reasonings show that the only possibility is $0\le t_{1}\le x\le t_{2}\le 1$, with $t=t_{2}$ maximum point on $[x, 1]$ for $H_{n, k}(t)$.

Similarly, for $t\in [0, x]$, since $H_{n, k}(t)=(x-t)t^{k}(1-t)^{n-k}$, using the above reasonings we obtain
$H^{\prime}_{n, k}(t)=t^{k-1}(1-t)^{n-k-1}[t^{2}(n+1)-t(nx+k+1)+kx]$ and that $t_{1}$ is a maximum point of $H_{n, k}(t)$ on $[0, x]$.

In conclusion, with $t_{1}, t_{2}$ given by (\ref{eq5}), we get
$$A_{n, k}(x)=\max\{(t_{2}-x)t_{2}^{k}(1-t_{2})^{n-k}, \, (x-t_{1})t_{1}^{k}(1-t_{1})^{n-k}\}.$$

Case (ii). Suppose first that $k=0$. By the calculation from the case (i), for $t\in [x, 1]$ we get $0=t_{1}\le x\le t_{2}=\frac{nx +1}{n+1}\le 1$, $H_{n, 0}(t)\ge 0$ and $H_{n, 0}(x)=H_{n, 0}(1)=0$, which by similar graphical reasonings leads to the fact that the maximum of $H_{n, 0}(t)$ on $[x, 1]$ is $H_{n, 0}(t_{2})=(t_{2}-x)(1-t_{2})^{n}$. Therefore, we recapture the case (i) with the convention that $0^{0}=1$. Similarly, for $t\in [0, x]$, we get that
the maximum of $H_{n, 0}(t)$ is $H_{n, 0}(t_{1})=(x-t_{1})(1-t_{1})^{n}$

The subcase $k=n$ is similar, which proves the lemma.
$\hfill \square$

{\bf Proof of Theorem 4.1.} According to Theorem 3.1, (i), we have to estimate
$$D_{n, \Gamma_{n}}(\varphi_{x})(x)=\sum_{k=0}^{n}p_{n, k}(x)\cdot \frac{(C)\int_{0}^{1}|t-x|t^{k}(1-t)^{n-k}d P_{\lambda_{n, k}}(t)}{(C)\int_{0}^{1}t^{k}(1-t)^{n-k}d P_{\lambda_{n, k}}(t)}.$$
First of all, by Definition 2.2, (ii), we get
$$(C)\int_{0}^{1}t^{k}(1-t)^{n-k}d P_{\lambda_{n, k}}(t)=\int_{0}^{+\infty}P_{\lambda_{n, k}}(\{t\in [0, 1] ;
t^{k}(1-t)^{n-k}\ge \beta\})d \beta$$
$$=\int_{0}^{1}P_{\lambda_{n, k}}(\{t\in [0, 1] ; t^{k}(1-t)^{n-k}\ge \beta\})d \beta$$
$$=\int_{0}^{1}\sup\{\lambda_{n, k}(s) ; s\in \{t\in [0, 1] ; t^{k}(1-t)^{n-k}\ge \beta\}\}d \beta$$
$$=\frac{1}{k^{k}n^{-n}(n-k)^{n-k}}\cdot
\int_{0}^{1}\sup\{s^{k}(1-s)^{n-k} ; s\in \{t\in [0, 1] ; t^{k}(1-t)^{n-k}\ge \beta\}\}d \beta.$$
For simplicity, denote $E_{n, k}=k^{k}n^{-n}(n-k)^{n-k}$, where again we take $0^{0}=1$. Since for $\beta > E_{n, k}$ we have
$\{t\in [0, 1] ; t^{k}(1-t)^{n-k}\ge \beta\}=\emptyset$ and since we can take $\sup\{s^{k}(1-s)^{n-k} ; s\in \emptyset\}=0$, it follows
$$(C)\int_{0}^{1}t^{k}(1-t)^{n-k}d P_{\lambda_{n, k}}(t)$$
$$=\frac{1}{E_{n, k}}\cdot \int_{0}^{E_{n, k}}\sup\{s^{k}(1-s)^{n-k} ; s\in \{t\in [0, 1] ; t^{k}(1-t)^{n-k}\ge \beta\}\}d \beta$$
\begin{eqnarray}\label{eq11}
=\frac{1}{E_{n, k}}\cdot \int_{0}^{E_{n, k}}E_{n, k}d \beta=E_{n, k}.
\end{eqnarray}
On the other hand, denoting $A_{n, k}(x)=\sup\{|t-x|t^{k}(1-t)^{n-k};t\in [0, 1]\}$, by Remark 2.3, (iii), (v) and by Lemma 4.2, for $t_{1}< t_{2}$ in (\ref{eq5}) we obtain
$$(C)\int_{0}^{1}|t-x|t^{k}(1-t)^{n-k}d P_{\lambda_{n, k}}(t)\le (C)\int_{0}^{1}A_{n, k}(x)d P_{\lambda_{n, k}}(t)$$
$$=A_{n, k}(x)(C)\int_{0}^{1}1 d P_{\lambda_{n, k}}(t) = \max\{(t_{2}-x)t_{2}^{k}(1-t_{2})^{n-k}, \, (x-t_{1})t_{1}^{k}(1-t_{1})^{n-k}\}$$
$$\le (t_{2}-x)t_{2}^{k}(1-t_{2})^{n-k} + (x-t_{1})t_{1}^{k}(1-t_{1})^{n-k}.$$
Since $\frac{t_{2}^{k}(1-t_{2})^{n-k}}{k^{k}n^{-n}(n-k)^{n-k}}\le 1$,
$\frac{t_{1}^{k}(1-t_{1})^{n-k}}{k^{k}n^{-n}(n-k)^{n-k}}\le 1$ and
by Lemma 4.2 we get
$$\frac{A_{n, k}(x)}{k^{k}n^{-n}(n-k)^{n-k}}\le (t_{2}-x)\cdot \frac{t_{2}^{k}(1-t_{2})^{n-k}}{k^{k}n^{-n}(n-k)^{n-k}} + (x-t_{1})\cdot \frac{t_{1}^{k}(1-t_{1})^{n-k}}{k^{k}n^{-n}(n-k)^{n-k}}$$
$$\le t_{2}-t_{1}=\frac{\sqrt{\Delta}}{n+1}\le \frac{\sqrt{(nx-k)^{2}+2x(n-k)+2k(1-x)+1}}{n}$$
$$\le \sqrt{(x-k/n)^{2}+2x/n+(2k/n)\cdot (1-x)/n +1/n^{2}}$$
$$\le |x-k/n|+\sqrt{2x}/\sqrt{n}+(\sqrt{2k}/\sqrt{n})\cdot \sqrt{(1-x)/n}+1/n,$$
this immediately implies
$$D_{n, \Gamma_{n}}(\varphi_{x})(x)\le \sum_{k=0}^{n}p_{n, k}(x)(|x-k/n|+\sqrt{2x}/\sqrt{n}+\sqrt{2k/n}\cdot \sqrt{(1-x)/n}+1/n)$$
$$\le \frac{\sqrt{x(1-x)}}{\sqrt{n}}+\frac{\sqrt{2 x}}{\sqrt{n}}+\frac{\sqrt{2}\sqrt{x(1-x)}}{\sqrt{n}}+\frac{1}{n}
=\frac{(1+\sqrt{2})\sqrt{x(1-x)}+\sqrt{2}\sqrt{x}}{\sqrt{n}}+\frac{1}{n}.$$
Above we have used the well-known estimate
$\sum_{k=0}^{n}p_{n, k}(x)|x-k/n|\le \frac{\sqrt{x(1-x)}}{\sqrt{n}}$ and the Cauchy-Schwarz inequality for Bernstein polynomials, $B_{n}(f)(x)\le \sqrt{B_{n}(f^{2})(x)}$, applied for $f(t)=\sqrt{t}$.

Finally, applying Theorem 3.1, (i), the proof of Theorem 4.1 follows. $\hfill\square$

{\bf Remark 4.3.} For $\mu=\sqrt{m}$ with $m$ denoting the Lebesgue measure on $[0, 1]$, another particular case of the Bernstein-Durrmeyer operators for which the quantitative estimates in Theorem 3.1 are applicable would be, for example, $$D_{n, \mu}(f)(x)=\sum_{k=0}^{n}p_{n, k}(x)\cdot \frac{(C)\int_{0}^{1}f(t)t^{k}(1-t)^{n-k}d \mu}{(C)\int_{0}^{1}t^{k}(1-t)^{n-k}d \mu}.$$
It is worth noting that the uniform convergence of this $D_{n, \mu}(f)$ to $f$, follows directly from the general Theorem 3.2 in \cite{Gal-Opris}. Also, it can be obtained by using the nonlinear Feller kind scheme expressed by Theorem 3.1 in \cite{Gal-Ann} (combined with Remark 3.2 there), since by direct calculation we can show  that $D_{n, \mu}(e_{1})$ converges uniformly to $e_{1}$ and $D_{n, \mu}((\cdot - x)^{2})$ converges uniformly to $0$, on $[0, 1]$.

{\bf Remark 4.4.} Since the Bernstein-Durrmeyer-Choquet operators in this paper can be defined with respect to a family of Borel or Choquet measures, combined  in various ways, this fact offers a very high flexibility and generality, allowing to construct operators having even better approximation properties.
A first example for this flexibility is shown by Theorem 3.4 and Remark 3.6.

For the second example, let us replace in formula (\ref{eq8}) the family $\Gamma_{n}$ of measures of possibilities $P_{\lambda_{n, k}}$, $k=0, ..., n$, by the family consisting in the Dirac measures $\delta_{k/n}$, $k=0, 1, ..., n-1$, (which are Borel measures and therefore with the corresponding Choquet integrals reducing to the classical ones) together with a monotone, submodular, strictly positive set function  $\mu$. Then, denoting by $B_{n}(f)(x)$ the classical Bernstein operators, for $D_{n, \Gamma_{n}}$ in (\ref{eq8}) we get
$$D_{n, \Gamma_{n}}(f)(x)-f(x)=\left [\sum_{k=0}^{n-1}p_{n, k}(x)f\left (\frac{k}{n}\right )+x^{n}\cdot \frac{(C)\int_{0}^{1}f(t)t^{n}d\mu(t)}{(C)\int_{0}^{1}t^{n}d \mu(t)}\right ] -f(x)$$
$$=B_{n}(f)(x)-f(x) +x^{n}\left [\frac{(C)\int_{0}^{1}f(t)t^{n} d \mu(t)}{(C)\int_{0}^{1}t^{n} d \mu(t)}-f(1)\right ].$$
Suppose now that $f\ge 0$ is strictly increasing and strictly convex on $[0, 1]$ and, for example, that $\mu(A)=\sqrt{m(A)}$ or $\mu(A)=\sin[m(A)]$, with $m$ the Lebesgue measure. The strict convexity implies $B_{n}(f)(x)-f(x)>0$ for all $x\in (0, 1)$ and the property of $f$ to be strictly increasing easily implies
$$\frac{(C)\int_{0}^{1}f(t)t^{n}d \mu(t)}{(C)\int_{0}^{1}t^{n}d \mu(t)}-f(1)<\frac{f(1)\cdot (C)\int_{0}^{1}t^{n}d \mu(t)}{(C)\int_{0}^{1}t^{n}d \mu(t)}-f(1)=0.$$
Therefore, in this case we get
$$|D_{n, \Gamma_{n}}(f)(x)-f(x)|<\max\left \{|B_{n}(f)(x)-f(x)|, x^{n}\left |\frac{(C)\int_{0}^{1}f(t)t^{n}d \mu(t)}{(C)\int_{0}^{1}t^{n}d \mu(t)}-f(1)\right |\right \},$$
i.e. for $x\in (0, 1)$, $D_{n, \Gamma_{n}}(f)(x)$ approximates better than $B_{n}(f)(x)$.

Here it is clear that $B_{n}(f)(x)$ can also be viewed as the Bernstein-Durrmeyer operators in the case when $\Gamma_{n}$ is composed by the Dirac measures $\delta_{k/n}$, $k=0, ..., n$, where we  note that although the Dirac measures are not strictly positive, however the Bernstein-Durrmeyer operators attached to them are well defined. This fact contrasts with the classical case in \cite{Berd1} when $\Gamma_{n}$ is composed by only one set function, independent of $n$, and when the strict positivity of the set function is necessarily for the convergence (see Theorem 1 in \cite{Berd1}). In other words, the strict positivity of the set functions in Theorem 3.1 is not always necessary.

For another example, let us consider the genuine Bernstein-Durrmeyer-Cho\-qu\-et operators given by
$$U_{n, \Gamma_{n}}(f)(x)=p_{n, 0}(x)\cdot \frac{(C)\int_{0}^{1}f(t)(1-t)^{n}d \nu_{n, 0}}{(C)\int_{0}^{1}(1-t)^{n}d \nu_{n, 0}}+ p_{n, n}(x)\cdot \frac{(C)\int_{0}^{1}f(t)t^{n}d \nu_{n, n}}{(C)\int_{0}^{1}t^{n}d \nu_{n, n}}$$
$$+\sum_{k=1}^{n-1}p_{n, k}(x)\cdot \frac{(C)\int_{0}^{1}f(t)p_{n-2,k-1}(t)d \mu_{n-2, k-1}(t)}{(C)\int_{0}^{1}p_{n-2,k-1}(t)d \mu_{n-2, k-1}(t)},$$
where $\Gamma_{n}=\{\nu_{n, 0}, \nu_{n, n}, \mu_{n-2, k-1}, k=1, ..., n-1\}$.

Let us denote by $G_{n}(f)(x)$, the classical genuine Bernstein-Durmeyer operator (see, e.g., \cite{Gonska}).
Choosing in $\Gamma_{n}$ the set functions $\mu_{n-2, k-1}, k=1, ..., n-1$ as the Lebesgue measure, $\nu_{n, 0}=\delta_{0}$ (as Dirac measure) and $\nu_{n, n}$ as a monotone, submodular and strictly positive set function, we immediately obtain
$$U_{n, \Gamma_{n}}(f)(x)-f(x)=G_{n}(f)(x)-f(x)+x^{n}\left [\frac{(C)\int_{0}^{1}f(t)t^{n} d \nu_{n, n}(t)}{(C)\int_{0}^{1}t^{n} d \nu_{n, n}(t)}-f(1)\right ].$$
Since the strict convexity of $f$ implies $G_{n}(f)(x)-f(x)>0$ for all $x\in (0, 1)$ (see, e.g., Lemma 2.1, (iv) in \cite{Gonska}), similar reasonings with those for the previous example show that if $f$ is strictly convex and strictly increasing on $[0, 1]$ (and, for example, $\nu_{n, n}(A)=\sqrt{m(A)}$ or $\nu_{n, n}(A)=\sin[m(A)]$), then $U_{n, \Gamma_{n}}(f)(x)$ approximates better $f$ on on $(0, 1)$ than the classical genuine operator, $G_{n}(f)(x)$.

{\bf Remark 4.5.} Recall that in \cite{Gal-Ann}, Example 4.2, for the nonlinear Picard-Choquet operators we have obtained a general estimate similar to that for the classical Picard operators, while for particular functions of the form $f(x)=Me^{-Ax}$, $M, A >0$, we got there essentially better error estimates.

{\bf Remark 4.6.} In \cite{Li} applications of the classical Bernstein-Durrmeyer operators to learning theory are presented. Taking onto account the very recent applications of the Choquet integral to learning theory (see, e.g., \cite{Hull1} and the references therein), it becomes of interest to see for potential applications of the Bernstein-Durrmeyer-Choquet operators to learning theory. Also, taking into account the applications of the classical Bernstein-Durrmeyer operators in regression estimation in, e.g., \cite{Raf} and the very recent applications of the Choquet integral to regression model, see, e.g., \cite{Grab}, it would be of interest to see for possible applications of the Bernstein-Durrmeyer-Choquet operators to the regression model.

\end{document}